%Authors: J. Bastero, J. Bernues, and A. Pena

%Title: An extension of Milman's reverse Brunn-Minkowski inequality

%Filename: basterobernuespena.tex
%TeX: Plain
%Length: 28799
%Received Date: 1/4/95
%SubjectClass: 26D20
%Abstract: The classical Brunn-Minkowski inequality
%states that for  $A_1,A_2\subset\R^n$ compact,
%$$ |A_1+A_2|^{1/n}\ge
%|A_1|^{1/n}+|A_2|^{1/n}\eqno(1) $$ where
%$|\cdot|$ denotes the Lebesgue measure on
%$\R^n$.
%In 1986 V. Milman {\bf [Mil 1]} discovered
%that if $B_1$ and $B_2$ are
% balls there is always a relative position of
%$B_1$ and $B_2$ for which a perturbed  inverse
%of $(1)$ holds. More precisely:\lq\lq{\sl
%There exists a constant $C>0$
% such that for all $n\in\N$ and any balls
%$B_1,B_2\subset\R^n$ we can find a linear
%transformation $u\colon\R^n\to\R^n$ with
%$|{\rm det}(u)|=1$ and
%
%$$|u(B_1)+B_2|^{1/n}\le
%C(|B_1|^{1/n}+|B_2|^{1/n})"$$}
%
%The aim of this paper is to extend this
%Milman's result to a larger class of sets.
%

%Citation: Geom. Funct. Anal. 5 (1995), no. 3, 572--581.

%32   space        33 ! exclam. pt.   34 " double quote  35 # sharp
%36 $ dollar       37 % percent       38 & ampersand     39 ' prime
%40 ( left paren.  41 ) rt. paren.    42 * asterisk      43 + plus
%44 , comma        45 - minus         46 . period        47 / division
%58 : colon        59 ; semi-colon    60 < less than     61 = equal
%62 > greater than 63 ? question mark 64 @ at
%91 [ left bracket 92 \ backslash     93 ] right bracket 94 ^ caret
% 95 _ underline    96 ` left single quote
%123 { left brace  124 | vertical bar 125 } right brace  126 ~ tilda

\hsize=15.5true cm \vsize=23true cm \hfuzz=1pt
\font\fp=cmr8 \parskip=3pt plus 1pt minus 1pt

\def\R{\hbox{$I\kern-3.5pt R$}}
\def\N{\hbox{$I\kern-3.5pt N$}}
\def\E{\hbox{$I\kern-3.5pt E$}}
\def\P{\hbox{$I\kern-3.5pt P$}}   \def\e{\varepsilon}
\def\b{\bigskip} \def\m{\medskip}
\def\d{\displaystyle }

\def\Pf{\noindent{\it Proof: }} \def\fin{
\hfill /// \vskip 12pt}

 \def\s{\sum_{i=1}^n}

\def\pconv{\hbox{$p$-{\rm conv}\ }}

\def\norm{\Vert\cdot\Vert}

\centerline {\bf An extension of Milman's
reverse Brunn-Minkowski inequality}

\vskip 6 pt \centerline {by} \vskip 6pt
\centerline {Jes\'us Bastero \footnote 
*{Partially supported   by Grant DGICYT PS
90-0120}, Julio Bernu\'es * and Ana Pe\~na *}
\footnote{}{AMS Class.:  46B20, . Key words:
$p$-convexity, entropy} \vskip 6pt
\centerline{\fp Departamento de Matem\'aticas.
Facultad de Ciencias} \centerline{\fp
Universidad de Zaragoza} \centerline{\fp
50009-Zaragoza (Spain)}  \vskip 20pt

\bigskip

\beginsection { 0. Introduction}

The classical Brunn-Minkowski inequality
states that for  $A_1,A_2\subset\R^n$ compact,
$$ |A_1+A_2|^{1/n}\ge
|A_1|^{1/n}+|A_2|^{1/n}\eqno(1) $$ where
$|\cdot|$ denotes the Lebesgue measure on
$\R^n$. 
 Brunn {\bf [Br]} gave the first proof of this
inequality for $A_1, A_2$  compact  convex
sets, followed by an analytical proof by
Minkowski {\bf [Min]}. The  inequality $(1)$
for compact sets, not necessarily convex,  was
first proved by Lusternik {\bf [Lu]}. A very
simple proof of it can be found in  {\bf [Pi
1]}, Ch. 1. 

It is easy to see that one cannot expect the
reverse inequality to hold
  at all, even if it is perturbed by a fixed
constant and we restrict ourselves to balls 
(i.e. convex symmetric compact sets with the
origin as an interior point).  Take for
instance $A_1=\{(x_1\dots x_n)\in\R^n\mid
|x_1|\le\e, |x_i|\le 1,2\le  i\le n\}$ and
$A_2=\{(x_1\dots x_n)\in\R^n\mid |x_n|\le\e,
|x_i|\le 1,1\le  i\le n-1\}$.

In 1986 V. Milman {\bf [Mil 1]} discovered
that if $B_1$ and $B_2$ are
 balls there is always a relative position of
$B_1$ and $B_2$ for which a perturbed  inverse
of $(1)$ holds. More precisely:\lq\lq{\sl
There exists a constant $C>0$
 such that for all $n\in\N$ and any balls
$B_1,B_2\subset\R^n$ we can find a linear
transformation $u\colon\R^n\to\R^n$ with
$|{\rm det}(u)|=1$ and 

$$|u(B_1)+B_2|^{1/n}\le
C(|B_1|^{1/n}+|B_2|^{1/n})"$$}

The nature of this reverse Brunn-Minkowski
 inequality is absolutely different from
others (say reverse Blaschke-Santal\'o
inequality, etc.).  Brunn-Minkowski
inequality   is an isoperimetric inequality,
 (in $\R^n$ it is its first and most important
consequence till now) and there is no inverse
to isoperimetric inequalities. So, it was a new
 idea that in the class of affine images of
convex bodies there is some kind of inverse.

 The result proved by Milman used   hard
technical tools (see {\bf [Mil 1]}).
 Pisier in {\bf [Pi 2]} gave  a new proof by
using interpolation and entropy
 estimates. Milman in {\bf [Mil 2]} gave
another proof by using the \lq\lq convex
surgery" and
 achieving also some entropy estimates.

The aim of this paper is to extend this
Milman's result to a larger class of sets. 
Note that simple examples show that some
conditions on a class of sets are clearly
necessary.

For $B\subset\R^n$   body (i.e. compact, with
non empty interior), consider
 $B_1=B-x_0$, where $x_0$ is an interior
point. If we denote by $N(B_1)=\cap_{\vert
a\vert\geq 1}aB_1$ the  balanced kernel of
$B_1$, it is clear that $N(B_1)$ is a balanced
 compact neighbourhood of the origin, so 
there exists $c>0$ such that  $B_1+B_1\subset
cN(B_1)$.  The Aoki-Rolewicz theorem (see 
{\bf [Ro], [K-P-R]}) implies that
  there is $0<p\le 1$, namely
$p=\log_2^{-1}(c)$, such that $B_1\subset \bar
B\subset  2^{1/p} B_1$, where $\bar B$ is  the
unit ball of some $p$-norm. This observation
will allow us to work in a $p$-convex
enviroment.

 The above construction allows us to define
the following parameter. For $B$ a body let
$p(B)$, $0<p(B)\leq 1$, be the supremum of
the  $p$ for which there exist a measure
preserving affine transformation of $B$,
$T(B)$, and a $p$-norm with unit ball $\bar B$
verifying $T(B)\subset \bar B$ and $\vert \bar
B\vert\leq \vert 8^{1/p} B\vert$, (by suitably
adapting
the results appearing in {\bf [Mil 2]}, it is
clear that $p(B)\geq p$ for any $p$-convex body $B$). 

Our main theorem is, 

\proclaim Theorem 1. Let $0<p\leq 1$. There
exists  $C=C(p)\geq 1$ such that for all
$n\in\N$ and all $A_1,A_2\subset\R^n$  bodies
such that $p(A_1),  p(A_2)\geq p$, there
exists  an affine transformation
$T(x)=u(x)+x_0$ with $x_0\in\R^n$,
$u\colon\R^n\to\R^n$ linear and $|{\rm
det}(u)|=1$ such that  $$|T(A_1)+A_2|^{1/n}\le
C(|A_1|^{1/n}+|A_2|^{1/n})$$

In particular, for the class of $p$-balls the
constant $C$ is universal (depending only on
$p$).

We prove this theorem in section 2. The key is
to estimate certain entropy numbers. We will
use the convexity of quasi-normed spaces of
Rademacher type $r>1$,  as well as
interpolation results
 and iteration procedures. 

We want to thank  Gilles Pisier for a useful
conversation
 during the preparation of this paper.

\b

\beginsection {1. Notation and background}

Throughout the paper $X, Y, Z$ will denote
finite dimensional real vector spaces.
 A quasi-norm on a real vector space $X$ is a
map $\Vert\cdot\Vert\colon X\to\R^{+}$  such
that\m

\item{i)} $\Vert x\Vert>0$ $\forall\ x\ne 0$.
\item{ii)} $\Vert tx\Vert=|t|\ \Vert x\Vert$
$\forall\ t\in\R, x\in X$. \item{iii)}
$\exists C\ge 1$ such that $\Vert x+y\Vert\leq
C(\Vert x\Vert+\Vert y\Vert)$ $\forall\ x,y\in
X$\smallskip

  If iii) is substituted by \smallskip

\item{iii')}$\Vert x+y\Vert^p\leq \Vert
x\Vert^p+\Vert y\Vert^p$ for  $x,y\in X$ and
some $0<p\le 1,$\smallskip

\noindent $\Vert\cdot\Vert$ is called a
$p$-norm on $X$. Denote by
 $B_X$ the unit ball of a quasi-normed or a
$p$-normed space. 

The above observations concerning the
$p$-convexification of our  problem can be
restated using $p$-norm and quasi-norm
notation. Recall that any compact balanced set
with 0 in its interior is the unit ball of a
quasi-norm. 

 By the concavity of the function $t^p$, any
$p$-norm is a quasi-norm with $C=2^{1/p-1}$.
Conversely, by the  Aoki-Rolewicz theorem, for
any quasi-norm with constant $C$ there exists
$p$, namely $p= \log_2^{-1}(2C)$, and a
$p$-norm $|\cdot|$ such that $|x|\le\Vert
x\Vert\le 4^{1/p}|x|,\
 \forall\ x\in X$. 

 A set $K\subset X$ is called $p$-convex if
$\lambda x+\mu y$, whenever $x,y\in K$, $
\lambda, \mu\ge 0,\ \lambda^p+\mu^p=1$. Given
$K\subseteq
 X$, the $p$-convex hull (or $p$-convex
envelope) of $K$ is the intersection
 of all $p$-convex sets that contain $K$. It
is denoted by $\pconv(K)$.  The closed unit
ball of a $p$-normed space
$(X,\Vert\cdot\Vert)$ will simply be
 called a $p$-ball. Any symmetric compact
$p$-convex set in $X$ with the
 origin as an interior point is the $p$-ball
associated to some $p$-norm.

We say that a quasi-normed space
$(X,\Vert\cdot\Vert)$ is of  (Rademacher) type
$q, 0<q\le 2$ if for some constant $T_q(X)>0$
we have

$${1\over 2^n}\sum_{\varepsilon_i=\pm
1}\Vert\s\varepsilon_i x_i \Vert\le T_q(X)
(\s\Vert x_i\Vert^q)^{1/q},\ \ \forall x_i\in
X, 1\le i\le n,\
 \forall n\in\N$$

 Kalton, {\bf [Ka]}, proved that any
quasi-normed space $(X, \Vert\cdot\Vert)$ of
type $q>1$ is convex. That is, the quasi-norm 
$\norm$ is equivalent to a norm and moreover,
the equivalence constant depends
 {\sl only}  on $T_q(X)$, (for a more precise
statement and proof of this fact see [{\bf
K-S}]). 

 Given $f,g\colon\N\to\R^+$ we write $f\sim g$
if there
 exists a constant $C>0$ such that
$C^{-1}f(n)\le g(n)\le Cf(n),\,  \forall\
n\in\N$. Numerical constants will always be
denoted by $C$ (or $C_p$ if
 it  depends only on $p$) although their value
may change from line to line.

Let $u\colon X\to Y$ be a linear map between
two  quasi-normed spaces and $k\ge 1$. Recall
the definition of the following numbers:
\smallskip

\item{}{\sl Kolmogorov numbers}: $\d d_k(u)=
\inf\{\Vert Q_S\circ u\Vert\mid S\subset Y \ \
{\rm subspace\ and}\ \ {\rm dim}(S)<k\}$
 where $Q_S\colon Y\to Y/S$ is the quotient
map.  \smallskip

\item{}{\sl Covering numbers}: For
$A_1,A_2\subset X$,  $N(A_1,A_2)=\inf\{N\in\N
\mid \exists\ x_1\dots x_N\in X\ {\rm such\
that}\  A_1\subset
 \bigcup_{1\le i\le N}(x_i+A_2)\}$.  \smallskip

\item{}{\sl Entropy numbers}: $\d
e_k(u)=\inf\{\e>0\mid N(u(B_X), \varepsilon
B_Y)\le 2^{k-1}\}$ \smallskip

The sequences   $\{d_k(u)\}, \{e_k(u)\}$ are 
non-increasing and satisfy
$d_1(u)=e_1(u)=\Vert u\Vert$. If ${\rm
dim}(X)={\rm dim}(Y)=n$ then $d_k(u)=0$ for
all $k>n$. Denote $s_k$ either  $ d_k$ or 
$e_k$. For all linear operators $u\colon X\to
Y$, $v\colon Y\to Z$
 we have  $\d s_k(v\circ u)=\Vert u\Vert
s_k(v)$ and $\d s_k(v\circ u)=\Vert  v\Vert
s_k(u), \ \forall k\in \N$ (called the ideal
property of $s_k$) and  $$s_{k+n-1}(v\circ
u)\le s_k(v) s_n(u) \ \ \ \forall\ k, n\in \N$$

The following two lemmas contain useful
information about these
 numbers. The first one extends to the
$p$-convex case its
 convex analogue due to Carl ({\bf [Ca]}). Its
proof mimics the ones of Theorem 5.1 and 5.2
in {\bf [Pi 1] } (see also [{\bf T}]) with minor changes.
In particular we identify $X$ as a quotient of $\ell_p(I)$,
for some $I$, and apply the metric lifting property of
$\ell_p(I)$ in the class of $p$-normed spaces (see Proposition
C.3.6 in [{\bf Pie]}).
 The second one contains easy facts about 
$N(A,B)$ and its  proof is similar to  the one
of Lemma 7.5. in  {\bf [Pi 1] }.

\proclaim Lemma 1. For all $\alpha>0$ and
$0<p<1$
 there exists a constant $C_{\alpha,p}>0$ such
that  for all linear map $u\colon X\to Y$,
$X,Y$  $p$-normed spaces and for all $n\in \N$
we have  $$\sup_{k\le n} k^{\alpha}e_k(u)\le
C_{\alpha,p} \sup_{k\le n} k^{\alpha}d_k(u) $$

\proclaim Lemma 2.  \item{i)} For all $A_1,
A_2, A_3\subset X$, $N(A_1,A_3) \le N(A_1,A_2)
N(A_2,A_3)$ \item{ii)} For all $t>0$ and
$0<p<1$ there is
 $C_{p,t}>0$ such that for all $X$ $p$-normed
space  of dimension $n$, $N(B_X,tB_X)\le
C^n_{p,t}$.  \item{iii)} For any 
$A_1,A_2,K\subset \R^n$,  $|A_1+K|\le
N(A_1,A_2) |A_2+K|$.  \item{iv)} Let $B_1,B_2$
be $p$-balls in $\R^n$
 for some $p$ and $B_2\subset B_1$; then
$\d{|B_1|\over |B_2|}\sim N(B_1,B_2).$

For any $B\subseteq\R^n$ $p$-ball the polar
set of $B$ is defined as
$$B^{\circ}:=\{x\in\R^n\mid \langle
x,y\rangle\le1,\,\, \forall\, y\in B\}$$ where
$\langle\cdot,\cdot\rangle$ denotes the
standard scalar  product on $\R^n$. Given
$B,D$ $p$-balls in $\R^n$ we define the
following  two numbers: $$s(B):=(|B|\cdot
|B^{\circ}|)^{1/n}$$ and
$$M(B,D):=\left({|B+D|\over |B\cap
D|}\cdot{|B^{\circ} +D^{\circ}|\over
|B^{\circ}\cap D^{\circ}|}\right)^{1/n}$$ \m

Observe that for any linear isomorphism
$u\colon\R^n \to\R^n$ we have $s(u(B))=s(B)$
and $$M(u(B),u(D)) =M(B,D).$$
 Recall that  $s( B_{\ell_p^n})\sim
n^{-1/p}\sim s( B_{\ell_2^n})^{1/p}, 0<p\le 1$
({\bf [Pi 1]} pg. 11). \m

The following estimates on these numbers are
known:\m

\item{a)} {\bf [Sa]}. For every symmetric
convex body  $B\subset\R^n$, $\d s(B)\le
s(B_{\ell_2^n})$ with equality  only if $B$
is  an ellipsoid. ({\sl Blaschke-Santal\'o's
inequality}). \m \item{b)} {\bf [B-M]}. There
exists a numerical constant
 $C>0$ such that for any $n\in\N$ and any   
symmetric convex body $B\subset\R^n$, $s(B)\ge
C s(B_{\ell_2^n})$.   \m \item{c)} {\bf [Mil
1]}. There exists a numerical constant
 $C>0$ such that for any $n\in\N$ and any  
symmetric convex body $B\subset\R^n$, there is
an ellipsoid (called Milman ellipsoid)
$D\subset\R^n$ such that $M(B,D)\le C$,({\sl Milman
ellipsoid theorem}).

 \b \beginsection {2. Entropy estimates and
reverse Brunn-Minkowski inequality}

We first introduce some useful notation: Let
$B_1, B_2\subset\R^n$ be two $p$-balls and
$u\colon\R^n\to\R^n$ a linear map. 
 We denote $u\colon B_1\to B_2$ the operator
between $p$-normed spaces  $u\colon(\R^n,
\norm_{B_1})\to(\R^n, \norm_{B_2})$ where
$\norm_{B_i}$ is the
 $p$-norm on $\R^n$ whose unit ball is $B_i$.

\m

{\bf Proof of Theorem 1:} \m

Let  $A_1, A_2$ be two bodies in $\R^n$ such
that $p(A_1)$, $p(A_2) \geq p$. It's clear
from the definition that there exist two
$\bar p$-balls, $B_1$, $B_2$, (for instance, 
$\bar p=p/2$) and two
measure preserving affine transformations $T_1$,
$T_2$, verifying  $$ \vert
T_2^{-1}T_1(A_1)+A_2\vert\leq \vert
B_1+B_2\vert $$ and $$ \vert
B_1\vert^{1/n}+\vert B_2\vert^{1/n} \leq C_p
\left(\vert A_1\vert^{1/n} +\vert
A_2\vert^{1/n}\right). $$ So, we only have to
prove the theorem for $p$-balls.
 
In the convex case a way to obtain the reverse
Brunn-Minkowski inequality is to prove that, for any
symmetric convex body $B$, there exists an ellipsoid $D$
verifying $\vert B\vert=\vert D\vert$ and 
$$\vert B+\Delta\vert^{1/n}\leq C
\vert D+\Delta\vert^{1/n}\eqno (2)$$
for any, say compact, subset $\Delta\subseteq\R^n$
($C$ is an universal constant independent of $B$ and $n$).

Indeed, let $B_1,B_2$ be two balls in $\R^n$. Suppose
w.l.o.g. that $D_i$, the ellipsoids associated to $B_i$
satisfy $u_2D_i=\alpha_iB_{\ell^n_2}$, where $u_i$
are linear mappings with $\vert {\rm det\ }u_i\vert=1$ and
$\vert
B_i\vert^{1/n}=\alpha_i\vert B_{\ell^n_2}\vert^{1/n}$.Then 
$$\eqalign{
\vert u_1B_1+u_2B_2\vert  ^{1/n}& \leq
C^2\vert u_1D_1+u_2D_2\vert^{1/n}  \cr
& = C^2(\alpha_1+\alpha_2)\vert B_{\ell^n_2}\vert^{1/n} =
C^2(\vert B_1\vert^{1/n}+\vert B_2\vert^{1/n})
}$$

 In view of the  preceding comments  and of straightforward
computations deduced from Lemma 2, in order to obtain $(2)$
for $p$-balls it is sufficent
 to associate an ellipsoid $D$ to each $p$-ball
$B\subset\R^n$ in such a way that the corresponding covering
numbers verify $N(B,D), N(D,B)\le C^n$ for
some constant $C$ depending only on $p$. 

It is important to remark now the fact that, what we deduce
from covering numbers estimate is that  the ellipsoid $D$
associated to $B$ actually verifies the stronger assertion
$$ C^{-1}\vert B+\Delta\vert^{1/n}\leq 
\vert D+\Delta\vert^{1/n}\leq C\vert B+\Delta\vert^{1/n}$$
for any compact set $\Delta$ in $\R^n$, with constant
depending only on $p$. Furthermore, the role  of  the
ellipsoid can be played by any fixed $p$-ball in a
\lq\lq spetial  position".

 Denote by $\hat B$ the convex hull of $B$. 

 By definition of $e_n$, if $e_n(id\colon B\to
D)\le \lambda$  then $N(B, 2\lambda D)\le
2^{n-1}$ and by Lemma 2-ii), $N(B,D)\le
\lambda^n$.
 (Of course, the same can be done with
$N(D,B)$). Therefore our problem  reduces to 
estimating entropy numbers. What we are going
to prove is really a stronger result than we
need, in the line of  Theorem 7.13 of [{\bf Pi
2}].

\proclaim{ Lemma 3}.  Given $\alpha >
1/p-1/2$, there exists a constant $C= C(\alpha
, p)$ such that, for any $n\in \N$ and for any
$p$-ball $B \in \R^n$ we can find an ellipsoid
$D\in \R^n$ such that  $$ d_k(D\to B)+
e_k(B\to D)\leq C \left( {n\over
k}\right)^\alpha 
$$ 
for every $1\leq k\leq n$.

\noindent {\it Proof of the Lemma}.
 From Theorem 7.13 of {\bf [Pi 2]} we can
easily deduce the following fact:  There
exists a constant $C(\alpha)>0$
 such that for any $1\leq k\leq n$, $n\in\N$
and  any ball $\hat B\subset\R^n$,
 there is ellipsoid $D_0\subset\R^n$ such that
the identity operator
 $ id\colon\R^n\to\R^n$ verifies  
 $$d_k( id\colon D_0 \to \hat B)\le
C(\alpha)\Big({n\over k} \Big)^\alpha\qquad 
{\rm and}\qquad e_k( id\colon \hat B\to D_0
)\le C(\alpha)\Big( {n\over k}\Big)^\alpha\,
\eqno (3)$$

For simplicity, since we are always going to
deal with   the identity operator, we will
denote $ id\colon B_1\to B_2$ by $B_1\to B_2$.

  Let $D_0$ be the ellipsoid associated to
$\hat B$ in $(3)$. It is well known  {\bf
[Pe]}, {\bf [G-K]} that
 $B\subseteq \hat B\subseteq n^{1/p-1}B$. This
means $\Vert B\to \hat
 B\Vert\le 1$ and $\Vert \hat B\to  B\Vert\le
n^{1/p-1}$. Now, $(3)$ and  the ideal property
of $d_k$ and $e_k$ imply 
 $$d_k(D_0\to B)\le C(\alpha) n^{1/p-1}
\left({n\over k}\right)^\alpha\qquad {\rm
and}\qquad  e_k(B\to D_0)\le C(\alpha)\left( 
{n\over k}\right)^\alpha\qquad \forall\ k\le
n$$

This let us to introduce the  constant $C_n$
as the infimum of the constants $C>0$ for
which the conclusion of lemma 3 is true for
all $p$-ball in $\R^n$.
 Trivially  $C_n\leq C(\alpha )\left(1+n^{1/p-1}\right)$. Let
$D_1$ be an almost optimal ellipsoid such that  
$$\eqalign{ d_k(D_1\to B)&\leq 2C_n \left(
{n\over k}\right)^\alpha \cr e_k(B\to D_1)&\leq
2C_n \left( {n\over k}\right)^\alpha }\eqno (4)$$
 for
every $1\leq k\leq n$.

Use the real interpolation method with
parameters  $\theta, 2$ to interpolate the
couple $ id\colon B\to B$ and $ id\colon
D_1\to B$.  It is straightforward from its
definition that for $\d B_\theta\colon=(B,D_1)
_{\theta,2}$, we have
 $$ d_k(B_\theta\to B)\le \Vert B\to
B\Vert^{1-\theta}
 (d_k(D_1\to B))^{\theta}\qquad \forall\ k\le
n$$  and therefore,   $$ d_k(B_\theta\to B)\le
 \Big(2C_n\left({n\over
k}\right)^\alpha\Big)^{\theta}\qquad \forall\
k\le n$$

Write $\d\lambda=4C_n \left({n\over k}\right)^\alpha$. By
definition of the entropy numbers,
 there exist $x_i\in \R^n$ such that $\d
B\subset  \bigcup_{i=1}^{2^{k-1}}x_i+2\lambda
D_1$. But by perturbing $\lambda$ with an
absolute
 constant we can suppose w.l.o.g. that $x_i\in
B$. For all $z\in B$, there
 exists $x_i\in B$ such that $\Vert
z-x_i\Vert_{D_0}\le 2\lambda$. Also by 
$p$-convexity, $\Vert z-x_i\Vert_B\le
2^{1/p}$. 

 A general result (see {\bf [B-L]} Ch. 3.)
 assures the existence of a constant $C_p>0$
such that 
 $$\Vert x\Vert_{B_\theta}\le C_p\Vert
x\Vert_B^{1-\theta}\Vert
x\Vert_{D_1}^{\theta}.$$  Therefore, for all
$z\in B$, there exists
 $x_i\in B$ such
 that $\Vert z-x_i\Vert_{B_\theta}\le
C_p\lambda^{\theta}$ which means  $$e_k(B\to
B_\theta)\le C_p \Big(2C_n\left({n\over
k}\right)^\alpha\Big)^{\theta}.$$

Since $\alpha>1/p-1/2$, then  we can pick  $\d
\theta\in (0,1)$ such that $\d
{2(1-p)\over2-p}<\theta<\min\{1,1-1/2\alpha\}$.
Then   $\d B_\theta$ has
 Rademacher type strictly bigger than $1$
because $ {1-\theta\over p}+ {\theta\over 
2}<1$.

By Kalton's result quoted before, we can
suppose
 that $B_\theta$ is a ball and therefore we
can apply to it $(3)$ for $\gamma=\alpha
(1-\theta)>1/2$
 and assure the existence of another ellipsoid
$D_2$ such that 
 $$d_k(D_2\to B_\theta)\le C(\gamma)\left(
{n\over k}\right)^\gamma \qquad {\rm and}\qquad
  e_k(B_\theta\to
 D_2)\le C(\gamma)\left( {n\over
k}\right)^\gamma \qquad {\rm and}\qquad
\forall\ k\le n$$ 

Recall that  $d_{2k-1}(D_2\to B)\le 
 d_k(D_2\to B_\theta)d_k(B_\theta\to B)$ and
the same for the $e_k$'s. Thanks to
 the  monotonicity of the numbers $s_k$ we can
use the what is known about $s_{2k-1}$ for all
$s_k$. Using the estimates
 obtained above we get  $\forall\ k\le n$, 
$$d_k(D_2\to B)\le C(p,\alpha)2^\theta C_n^\theta 
\left({n\over
k}\right)^{\gamma+\alpha\theta}\qquad  {\rm
and}\qquad   e_k(B\to D_2)\le C(p,\alpha) 2^\theta
C_n^\theta  \left({n\over
k}\right)^{\gamma+\alpha\theta}. $$ Hence by
the election of $\gamma $ and by minimality we
obtain $ C_n^{1-\theta}\leq C(p,\alpha)2^\theta$, and the
conclusion of the lemma holds.

\fin

The theorem follows  now from the estimate we
achieved  in Lemma 3 and by Lemma 1. Indeed, given any
$\alpha>1/p-1/2$,
 if $D$ is the ellipsoid associated to $B$ by
Lemma 3, we
have   $$\eqalign { n^{\alpha}e_n(D\to B) &
\le \sup_{k\le n}k^{\alpha}e_k(D\to B) 
 \le C(\alpha, p)\sup_{k\le
n}k^{\alpha}d_k(D\to B)\le C(\alpha ,
p)\sup_{k\le n}k^{\alpha}{n^{\alpha}\over
k^{\alpha}}\cr &=C(\alpha , p)n^{\alpha} \cr
}$$  and so, $e_n(D\to B)\le C(\alpha,p)$.
On the other
 hand just take $k=n$ in Lemma 3 and so,
$e_n(B\to D)\le C(\alpha,p)$. 

 Finally
observe that since the constant $C(\alpha,p)$ 
depends only on $p$ and $\alpha$ and we can
take any $\alpha>1/p-1/2$ the thesis of the
theorem as stated inmediately follows.

\fin

\beginsection {3. Concluding remarks}

We conclude this note by stating the
corresponding versions of a)
Blaschke-Santal\'o, b) reverse
Blaschke-Santal\'o and c) Milman ellipsoid
theorem, cited in section 1, in the context
of $p$-normed spaces.

\proclaim Proposition 1. Let $0<p\leq 1$. 
There exists a numerical  constant $C_p>0$ 
such that for every
 $p$-ball $B\subseteq \R^n$, 
 $$ C_p \big(s(B_{\ell_2^n})\big)^{1/p}\leq
s(B)\leq s(B_{\ell_2^n})$$ and in the second
inequality, equality holds  if only if $B$ is
an ellipsoid

\Pf Denote by $\hat B$ the convex envelope of
$B$.  Since $\hat B^{\circ}=B^{\circ}$  we
have $s(B)\le s(\hat B)\le
 s(B_{\ell_2^n})$.  If  $\d s(B)=
s(B_{\ell_2^n})$, then $\hat B$ is an
 ellipsoid. We will show that  $B=\hat B$. 
Every $x$ in the boundary of $\hat B$
 can be written as $x=\sum \lambda_ix_i,
x_i\in B$, $\sum\lambda_i=1$; but  since $\hat
B$ is an ellipsoid, $x$ is an extreme point of
$\hat B$ and so $x=x_i$  for some $i$ that is
$x\in B$. This shows $B=\hat B$ and we are
done. 

For the first inequality, $B\subseteq \hat
B\subseteq n^{1/p-1}B$  easily implies 
 $\d \left(|\hat B|\over |B|\right)^{1/n} \le
n^{1/p-1}$ and so,   $$\eqalign{ s(B) &
=(|B|\cdot |B^{\circ}|)^{1/n}=(|B|\cdot | \hat
B^{\circ}|)^{1/n}= \left(|B|\over |\hat
B|\right)^{1/n}(|\hat B|\cdot 
|B^{\circ}|)^{1/n}\ge {s(\hat B)\over
n^{1/p-1}}\cr & \ge {C s(B_{\ell_2^n})\over
n^{1/p-1}}= C n^{-1/p}= C_p
\big(s(B_{\ell_2^n})\big)^{1/p}  \cr } $$

\fin

 The  left inequality above is sharp since
$s(B_{\ell_p^n})=C_p\big(s(B_{\ell_2^n})\big)^{1/p}$.
The  right
 inequality is also  sharp since every ball is
a $p$-ball for
 every $0<p<1$. And it is sharp even if we
restrict
 ourselves to the class of $p$-balls which are
not $q$-convex for any $q>p$,  as it is showed
by the following example:
  Let $\varepsilon>0$ and $C_\varepsilon$ be a
relatively open cap in $S^{n-1}$ centered in
$x=(0,\dots,0,1)$ of radius  $\varepsilon$.
Write $K=S^{n-1}\setminus
\{C_\varepsilon\cup-C_\varepsilon\}$. The
$p$-ball $\pconv(K)$ is not $q$-convex for any
$q>p$ and we can pick $\varepsilon$ such that
$\d{s(\pconv(K))\over s(B_{\ell_2^n})}\sim 1 $.

 Observe that the left inequality is actually
{\sl equivalent}
 to the existence of a constant $C_p>0$ such
that for every $p$-ball $B$, $\d \left(|\hat
B|\over |B|\right)^{1/n}\le C_p n^{1/p-1}$ 
and by Lemma 2 iv),  this is also equivalent
to the inequality $N(\hat B, B)\le C_p
n^{1/p-1}$.

With respect to the Milman ellipsoid theorem
we obtain

\proclaim Proposition 2.  Let $0<p<1$. There
exists
 a numerical constant $C_p>0$  such that for
every $p$-ball $B$ there is an
 ellipsoid $D$ such that $M(B,D)\le  C_p
n^{1/p-1}$. 

\Pf Given a $p$-ball $B$ let $D$ be the
Milman  ellipsoid of $\hat B$. Then, 
$$\eqalign{ M(B,D) & =\left({|B+D|\over |B\cap
D|}\cdot{|B^{\circ}+D|\over |B^{\circ}\cap
D|}\right)^{1/n}\cr &= \left({|\hat B+D|\over
|\hat B\cap D|}\cdot {|\hat B^{\circ}+D|\over
|\hat B^{\circ}\cap D|}\right)^{1/n}
\left({|B+D|\over |\hat
B+D|}\right)^{1/n}\left({|\hat B\cap D|\over
|B\cap D|}\right)^{1/n}\cr  &\le M(\hat B, D)
\left({|\hat B\cap D|\over |B\cap
D|}\right)^{1/n}\le C_p n^{1/p-1} \cr }$$
\hfill ///

 The bound for $M(B,D)$ is sharp. Indeed, if
there
 was a function $f(n)<<n^{1/p-1}$ such that
for every a $p$-ball $B$
 there was an ellipsoid $D$ with $M(B,D)\le
f(n)$, then  $${s(B_{\ell_2^n})\over
f(n)}={s(D)\over f(n)}\le \left(|B\cap
D|\cdot|B^{\circ}\cap D|\right)^{1/n}\le
s(B)$$ and we would have, $s(B)\ge
{s(B_{\ell_2^n})\over f(n)}>>n^{-1/p}$ which
is not possible.

{\bf Acknowledgments.} The authors are indebted to the
referee for some fruitful comments which led them to improve
the presentation of the paper.

\b

\beginsection References

\item{\bf[B-L]} Bergh, J. L\"ofstr\"om, J.:
{\sl
 Interpolation Spaces. An Introduction}.
Springer-Verlag, 223. Berlin Heildelberg 
(1976).  \m

\item{\bf[B-M]} Bourgain, J., Milman, V.D.:
{\sl On  Mahler's conjecture on the volume of
a convex symmetric body and its polar}.
 Preprint I.H.E.S., March 1985. \m

\item{\bf[Br]} Brunn, H.: {\sl \"Uber Ovale
und  Eifl\"achen}. Innaugural dissertation.
M\"unchen 1887.  \m

\item{\bf[Ca]} Carl B.: {\sl Entropy numbers,
$s$-numbers,  and eigenvalue problems}. J.
Funct. Anal. {\bf 41}, 290-306 (1981). \m

\item{\bf[G-K]}  Gordon, Y., Kalton, N.J.:
{\sl Local
 structure for quasi-normed spaces}. Preprint.
\m 

\item{\bf[Ja]} Jarchow, H.: {\sl Locally
convex spaces}. B.G. Teubner Stuttgart,
(1981). \m

\item{\bf[Ka]} Kalton, N.: {\sl Convexity,
Type and
 the Three Space Problem}.  Studia Math. {\bf
69}, 247-287 (1980-81). \m

\item{\bf[K-S]} Kalton, N., Sik-Chung Tam:
{\sl Factorization theorems for quasi-normed
spaces}. Houston J. Math. {\bf 19} (1993),
301-317. \m

\item{\bf[K-P-R]} Kalton, N., Peck, N.T.,
Roberts, J.W. : {\sl An F-space sampler}. 
London Math. Soc. Lecture Notes 89. Cambridge
Univ. Press. Cambridge (1985). \m

\item{\bf[Lu]} Lusternik, L.A.: {\sl Die
Brunn-Minkowskische  Ungleichung f\"ur
beliebige messbare Mengen}. C.R. Acad.Sci.URSS
{\bf 8}, 55-58 (1935). \m

\item{\bf[Mil 1]} Milman, V.D.: {\sl
In\'egalit\'e de Brunn-Minkowsky  inverse et
applications \'a la th\'eorie locale des
espaces norm\'es}. C.R. Acad.Sci.Paris {\bf
302},S\'er 1, 25-28 (1986). \m

\item{\bf[Mil 2]} Milman, V.D.: {\sl
Isomorphic symmetrizations and  geometric
inequalities}. In \lq\lq Geometric Aspects of
Functional
 Analysis" Israel Sem. GAFA. Lecture Notes in
Mathematics, Springer,
 {\bf 1317}, (1988), 107-131. \m

\item{\bf[Min]} Minkowski, H.:{\sl Geometrie
der Zahlen}.  Teubner, Leipzig (1910). \m

\item{\bf[Pe]} Peck, T.: {\sl Banach-Mazur
distances and projections on $p$-convex
spaces}. Math. Zeits. {\bf 177}, 132-141
(1981). \m

\item{\bf[Pie]} Pietsch, A.: {\sl Operator
ideals}. North-Holland, Berlin 1979.\m

\item{\bf[Pi 1]} Pisier, G.: {\sl The volume
of convex bodies and  Banach Space Geometry}.
Cambridge University Press, Cambridge (1989).
 \m

\item{\bf[Pi 2]} Pisier, G.: {\sl A new
approach to several results of V. Milman}. J.
reine angew. Math. {\bf 393}, 115-131 (1989).
\m

\item{\bf[Ro]} Rolewicz, S.: {\sl Metric
linear spaces}. MM 56.
 PWN, Warsaw (1972). \m

\item{\bf[Sa]} Santal\'o, L.A.: {\sl Un
invariante af\'{\i}n para los
 cuerpos convexos de $n$ dimensiones}.
Portugal Math. {\bf 8}, 155-161 (1949).
 \m \item{\bf[T]} Triebel, H.: {\sl Relations
between approximation numbers and entropy
numbers}. Journal of Approx. Theory {\bf 78},
(1994), 112-116.

\bye